\documentclass[12pt,a4paper]{article}

\usepackage[utf8]{inputenc}
\usepackage[english]{babel}
\usepackage{amsthm}
\usepackage{cite}
\usepackage{amsfonts}
\usepackage{amssymb}
\usepackage{amsmath}
\usepackage{graphicx}

\newcommand{\re}{\mbox{\rm Re}\,}
\newcommand{\im}{\mbox{\rm Im}\,}
\newcommand{\CC}{\mathbb{C}}
\newcommand{\const}{\mbox{const}}

\newenvironment{Definition}{\par\bigskip\noindent{\bf Definition.}}{\par\bigskip}

\newtheorem{Theorem}{Theorem}
\newtheorem{Lemma}{Lemma}

\title{The Limit Spectral Graph in the Semi-Classical Approximation for the
Sturm-Liouville Problem With a Complex Polynomial Potential}
\author{A.\,A.~Shkalikov,~S.\,N.~Tumanov}
\date{March 29th, 2016.}

\begin{document}

\maketitle
\begin{abstract}
The limit distribution of the discrete spectrum of the Sturm--Liouville problem
with complex--valued polynomial potential on an interval,
on a half--axis, and on the entire axis is studied. It is shown that at large
parameter values, the eigenvalues are concentrated
along the so--called limit spectral graph; the curves forming this graph are
classified. Asymptotics of eigenvalues
along curves of various types in the graph are calculated.
\end{abstract}
\frenchspacing

In this paper we study the limit distribution of the discrete spectrum in the
Sturm--Liouville problem when the parameter
$k\to+\infty$:
\begin{gather}
\label{Intr_Eq_1}
-y''(z)+k^2P(z,\lambda)y(z)=0, \\
\label{Intr_Eq_2} y(a) =y(b)=0.
\end{gather}

Here $\lambda$ is the generalized spectral parameter which  varies  in the  domain $G$
of the complex plane and $P$ is a polynomial
of degree $n\geqslant 1$ in  $z$ with coefficients  analytically depending on $\lambda \in G$:

$$
P(z,\lambda)=a_n(\lambda)z^n+a_{n-1}(\lambda)z^{n-1}
+\cdots +a_1(\lambda)z+a_0(\lambda).
$$

The  points  $a$ and $b$ at which the boundary conditions are defined  we assume
to be arbitrary complex numbers; one of them or both may be infinite. With the
infinite point $a$ we associate the ray $\Upsilon_\varphi = \{re^{i\varphi}, \
r>0\}$  and write $a= e^{i\varphi}\infty$, while the boundary condition
$y(a)=0$ in this case we understand in the following sense: $y= y(re^{i\varphi})\to 0$ as $r\to\infty$
along the ray $\Upsilon_\varphi$.
 The set of values
$\lambda\in G$ for which
 the equation \eqref{Intr_Eq_1}
has a  non-trivial solution satisfying the boundary conditions \eqref{Intr_Eq_2} we call the
spectrum of the problem under consideration. In the usual classical setting
 the points  $a$ and $b$  are assumed  to be  real or equal to $\pm\infty$
and in this case the problem is considered on the finite real segment $[a,b]$, on the
semi-axis $[a, \infty)$ or on the whole axis $(-\infty, \infty)$.
The spectral parameter is usually
regarded as linear (i.e. it is assumed that
$a_j(\lambda) = c_j =\const$, $j=1,\ldots,n$, $a_0(\lambda) = c_0-\lambda$).

We present here more general setting, however, it by no means  complicates the results obtained
below, and even makes them simpler formulated. It should be noted that when
talking about the semi-classical approximation one  typically uses a small
parameter $\varepsilon$ or $h$ and writes the equation in the following
form:
\[
-hy''+ (P(z) -\lambda)y =0,\ h\to 0.
\]
\noindent Surely, when $h=1/k^2 \to 0$, this is  equivalent to
\eqref{Intr_Eq_1} with $k\to\infty$.

In what follows, we exclude the
degenerate cases from consideration and require additionally for the
coefficient $a_n(\lambda)$ not to vanish in $G$, and for the
zeros
$z_j$, $j=1,\ldots,n,$ of the
 polynomial  $P$ to be
different functions of $\lambda$, or different germs of a function analytic
 in $G$ everywhere except some algebraic branching points.

This paper is a development of previous works of the authors \cite{Sh1, TSh1,
TSh2, TSh3, Sh2, PSh}.
References to the works of other authors related to the subject performed before 2004
are given in \cite{Sh2}. Among the more recent  works we mention, for example,
\cite{ESh, EG}. Here we
do not give more details since we plan to do it in a more detailed version of
this work.

In what follows we  assume that the reader is familiar with the method of phase integrals (WKB
method). For details we refer  the readers to the papers and books
\cite{EvgrafFedoruk,Fedoruk1,Fedoruk2,Fedoruk3,Heding}. Here we recall
the basic definitions and notations.

Let us consider the function
\[
S(z_0,z;\lambda)=\int\limits_{z_0}^z\sqrt{P(\zeta,\lambda)}\,d\zeta.
\]
Here the branch of a root is fixed, and the integration is carried out along
a path from $z_0$ to $z$, not crossing the turning points $z_j=
z_j(\lambda)$, that is, the  zeros of  $P(z)$.
The system of lines in $\mathbb C$, determined for each turning point (zero of $P$)
$z_j$:
\[
\Gamma_j=\Bigl\{z\in\CC\,\Bigl|\Bigr.\,\re S(z_j,z;\lambda)=0\Bigr\},
\]
is called {\it the Stokes complex} associated with the
 turning point $z_j$.
The maximal connected
components of $\Gamma_j$ starting from $z_j$ and containing no other turning
points are called {\it Stokes lines}. Each Stokes line is either infinite and
contains a turning point at the origin or is finite and connects two
different turning points. In the second case only the initial turning point
belongs  to the Stokes line. Thus, the Stokes lines are always understood as
oriented curves --- either infinite or finite curvilinear semi-intervals ---
starting at a turning point.
The assembly  $\mathfrak{G}(\lambda)$ of all the
Stokes complexes corresponding to all turning points is called {\it the Stokes
Graph}.
The  Stokes Graph splits the complex plane into simply connected infinite domains
called {\it basic domains}. Among them there are two types of domains. The
first type relates to domains that are of the {\it half--plane} type (the function $S$
 maps them conformally to a half--plane $\re S>\alpha$ or $\re
S<\alpha, \ \, \alpha\in\mathbb R$);  the second type relates to the so-called  {\it strip}
type domains (the function $S$ maps them to  vertical strips in $\mathbb C$).

A domain  $D$  of the complex plane is called {\it canonical} if the function $S$
 maps it one-to-one to the entire complex plane with a finite number of
vertical cuts. A canonical domain contains two half--plane type domains and
some (possibly empty) set of strip type domains (the common boundaries are also
included).

The interior of angles in $\mathbb C$ formed by the rays of the  form
\begin{gather}
 l_j=\Bigl\{ re^{i\varphi_j}-\frac{a_{n-1}}{na_n},\  r>0\Bigr\},\notag\\
  \label{Intr_Eq_lk}
   \varphi_j=\frac{\pi-\arg a_n}{n+2}+\frac{2\pi j}{n+2},\ \,  j=0,\ldots,n+1,
\end{gather}
are called {\it the Stokes sectors}. These and only these rays are the asymptotes
for infinite Stokes lines.

\begin{Definition}
\label{Def_Linked_1}
We say that a point $\lambda \in G$
is  {\it exceptional}  for the boundary
value problem \eqref{Intr_Eq_1}, \eqref{Intr_Eq_2},
if at least one of the boundary conditions is specified in an infinite
point, say, $a= e^{i\varphi}\infty$, and $\varphi$ coincides with one of
the values $\varphi_j$
 in \eqref{Intr_Eq_lk}, i.e. the direction of one of the infinite
points coincides with the direction of one of the asymptotes to the Stokes lines.
\end{Definition}

Further in the case $a_n(\lambda) = a_n = \text{const}$\ \, and $a= e^{i\varphi}\infty$  we always assume that
$ \varphi \ne \varphi_j$ where $\varphi_j$ are defined by \eqref{Intr_Eq_lk}. The same we assume for the value $\psi$
if the second boundary point $b=  e^{i\psi}\infty$ is also infinite.

\begin{Theorem}
\label{Def_Th_1}
 Let $g$ be a compact set in the domain $G$ such that all points
$\lambda \in g$  are not exceptional for the boundary value problem \eqref{Intr_Eq_1}, \eqref{Intr_Eq_2}.
Then the spectrum of this problem is discrete in a small neighborhood of $g$.
  \end{Theorem}

If $a_n(\lambda)\ne \text{const}$ and $a_n(\lambda) \ne 0$ then $u(\lambda):=\arg\, a_n(\lambda)$
is a nonconstant harmonic function on $G$.
Therefore, the conditions $\arg\, a_n(\lambda) = \pi-\varphi_j(n+2)$ determines smooth curves
in the domain $G$. In what follows, we eliminate the points of these curves
from consideration passing, is necessary, to a narrower
domain $G'\subset G$ containing  no exceptional  points.  Throughout the rest of the paper
we assume that any compact set in $G$ satisfies the condition of Theorem 1.

According to our assumption, the roots $z_j(\lambda)$ of the polynomial $P(z, \lambda)$
are different algebraic functions or germs of an analytic function with algebraic
branching points, which can accumulate  only to the boundary of $G$.
Further we wish to avoid the branching points and  pass to a  narrower domain $G'\subset G$
 which does not contain points of this type. Certainly, it may happen that the obtained domain $G'$
 is  multiply--connected. As before we will write $G$ instead of $G'$.

Our nearest goal is to define the curves in $G$ that play an important role
 in the description of the spectrum of the problem \eqref{Intr_Eq_1}, \eqref{Intr_Eq_2}
as $k\to\infty$.

Recall  that a Stokes complex is called {\it simple} if it is generated by a
simple turning point $z_j$ and all the three Stokes lines that go out from it
are infinite (i.e., this complex contains no other turning points). All the
other complexes are called {\it compound}. If in addition to $z_j$ a complex
contains $z_{j_2}$,\ldots,$z_{j_s}$ turning points  (with multiplicity taken into account),
then it is called the $s$-{\it
point} complex.
We say that a point $\lambda\in G$ is  {\it regular} if the Stokes
graph $\mathfrak{G}(\lambda)$
 consists of only simple complexes. All the other points we call {\it
singular}. If $\lambda$ is a singular point, then there exists a Stokes complex
which includes at least two different turning points, say, $z_j(\lambda)$
and $z_l(\lambda)$.
Let us consider the  integral
\begin{equation}
\label{Sjl}
 S_{j,l}(\lambda)
=\int\limits_{z_j(\lambda)}^{z_l(\lambda)}\sqrt{P(\zeta,\lambda)}\,d\zeta,
\end{equation}
where the integration is carried out along the Stokes line connecting $z_j$ and
$z_l$. It can be shown that $S_{j,l}(\lambda)$ admits the analytic continuation to
the domain $G$  (this continuation is multivalued and depends on the path in
$G$). Consequently, the   function $u(\lambda):= \re S_{j,l}(\lambda)$ is
harmonic in any simply connected domain belonging to $G$, and the  set
\begin{equation}
\label{gammajl}
 \gamma_{j,l} =\Bigl\{\lambda \in G\,\Bigl|\Bigr.\,\re S_{j,l}(\lambda)
=0\Bigr\}
\end{equation}
determines smooth curves in $G$, possibly with bifurcation points, when
$S'_{j,l}(\lambda) =0$.

It  follows from the representations
\eqref{Sjl} and \eqref{gammajl}  that for all points $\lambda \in
\gamma_{j,l}$  the Stokes graph contains a compound complex connecting the turning points
$z_j$ and $z_l$. We shall call the curves $\gamma_{j,l}$  {\it
singular curves}.
If different curves $\gamma_{j,l}$ and $\gamma_{j,s}$  intersect at a point
$\lambda_0$, then the Stokes graph $\mathfrak{G}(\lambda_0)$
 includes at least  three--point complex. Intersection of three singular curves at one
point generates a four-point complex, etc. Of course, for some indices $j,l$
the set $\gamma_{j,l}$ may be empty in $G$.
Thus, a set of singular points in $G$ consists of the curves $\gamma_{j,l}$ and
in a generic position the number of their intersection  points is finite in
any compact subset  $g\subset G$.
We do not exclude the situation when the intersection of different curves $\gamma_{j,l}$
and $\gamma_{p,s}$  forms a continuum.

Suppose that the point $a$  which defines the boundary condition is finite. We say that
a point  $\lambda \in G$  is {\it critical} with respect  to  $a$
 if there is a Stokes complex $\Gamma_j$
  with turning point $z_j(\lambda) \ne a$
such that  one of its Stokes lines intersects  the point  $a$. Let us
consider the following integral
\[
 C_{j}(a,\lambda) =\int\limits_{z_j(\lambda)}^a\sqrt{P(\zeta,\lambda)}\,d\zeta,
\]
where the integration is carried out along the  Stokes line connecting $z_j$  and
$a$. It can be shown that $C_{j}(a, \lambda)$ admits the analytic continuation in
$\lambda$ along any path in $G$ which does not intersects  the zeros of $P$. Therefore,
$u(\lambda): = \re C_{j}(a,\lambda)$ is a locally harmonic function and the set
\[
 \gamma_j(a) =\Bigl\{\lambda \in G\,\Bigl|\Bigr.\,\re C_{j}(a,\lambda) =0\Bigr\}
\]
determines the curves in $G$ which we shall call {\it critical curves}. If the
point $b$  is also finite, we define similarly the function $C_j(b, \lambda)$
and the critical curves  $\gamma_j(b)$
(for some indices the sets $\gamma_j(a)$ and
$\gamma_j(b)$ may be empty or coincide).

%In the case of finite $a$ and (or) $b$ we say that a point $\lambda \in G$ is
%{\it particular} if for some  $j$: $z_j(\lambda)=a$ or $z_j(\lambda)=b$ {\it particular
%points}. It can be shown that the $\mathcal{E}$ set of particular points is discrete.

Suppose now that both $a$ and $b$ boundary points are finite. A point  $\lambda\in
G$ will be called {\it balanced point} with respect to the boundary points  $a$ and $b$
 if both these points lie in a canonical domain determined by the Stokes graph
$\mathfrak G(\lambda)$  and
\begin{equation}\label{B}
\re B(a,b,\lambda) =0,\text{ where } B(a,b,\lambda) =
 \int\limits_a^b \sqrt{P(\zeta,\lambda)}\,d\zeta.
\end{equation}
Here the integration is carried out along the path which entirety lies in the
canonical domain. Evidently, the set of balanced points with respect
to $a$ and $b$ consists of smooth curves in $G$ which we shall call {\it
balanced curves}.
Let us denote this set by $\gamma (a,b)$.  We remark that the curves of this set take the origin
 from the intersection points of the critical curves $\gamma_j(a)$ and $\gamma_l(b)$
 (if there are such points in $G$), or appear at the boundary of $G$ and finish at
similar points of intersection of critical curves, or go to the boundary of  $G$
 (in particular, go to $\infty$ if $G=\mathbb C$).

\begin{Definition}
\label{Def_Linked_2}
Given  $\lambda \in G$ let $\mathfrak G(\lambda)$ be the Stokes graph
 corresponding to
this point, and  $\Gamma_j$  be a certain Stokes complex in this graph. The
domains which the complex $\Gamma_j$  splits the complex plane $\mathbb C$  into, we shall call
{\it basic domains} for this complex. A domain  $D\subset\mathbb C$
is said to be  {\it admissible} with respect to the complex $\Gamma_j$ if
$D$ intersects no more than two basic domains of this complex.
In particular, {\it maximal admissible domains} (domains that have no
nontrivial admissible extensions) consist of two neighboring basic domains of
the complex  $\Gamma_j$
 and the  Stokes line (excluding its starting point) which
is their common boundary.

We say that the boundary points $a$ and $b$
are {\it linked} with respect  to
the complex $\Gamma_j$  if they both lie in some admissible domain for $\Gamma_j$.
We imply that an infinite point $a= e^{i\varphi}\infty$
belongs to a domain if the ray $\arg z = \varphi$ is  asymptotically contained in
this domain.
\end{Definition}

The usefulness of the definition of linkedness with respect to complexes is
clarified by the following statement.
\begin{Lemma}
\label{VKB_Lm_CanonicD} Given  $\lambda \in G$ assume that the boundary
points  $a$ and $b$  do not coincide with the
turning points of the Stokes graph $\mathfrak G(\lambda)$. Then the points
$a$ and $b$  are
linked with respect to all the Stokes complexes $\Gamma_j\subset \mathfrak G(\lambda)$
if and only if both these points belong to a canonical
 domain determined by $\mathfrak G(\lambda)$.
\end{Lemma}

Balanced curves are determined by the position of the  boundary points $a$ and $b$
with respect to the  Stokes graphs $\mathfrak G(\lambda)$, critical curves depend only on one of the
boundary points and singular curves do not dependent on $a$ and $b$.  Next, we shall select
some parts of  critical and singular curves which play a crucial role
in the description of the spectrum of the problem in question as the parameter
$k\to\infty$.

In the case of a finite point $a$  we shall
denote by $\gamma_j^b(a)$ the set of points $\lambda \in \gamma_j(a)$
 for which the $a$ and $b$ are not linked
with respect to the Stokes complex $\Gamma_j\subset\mathfrak G(\lambda)$.
Similarly, in the case of a finite point $b$
we shall denote by $\gamma_j^a(b)$ the set of points $\lambda \in \gamma_j(b)$
 for which the points $b$ and $a$
  are not
linked with respect to the Stokes complex $\Gamma_j\subset\mathfrak G(\lambda)$.
Finally, by $\gamma_{j,l}^{a,b}$ we shall denote
the  set of points $\lambda\in\gamma_{j,l}$
 for which the $a$ and $b$  are not linked with respect to the Stokes complex
 $\Gamma_j= \Gamma_l\subset\mathfrak G(\lambda)$.
The selected parts  of critical and singular curves will be called
{\it essential} critical and singular ones.

\begin{Definition}
\label{Def_Linked_3}
A set $T \subset G$  is called limit spectral set of the  problem \eqref{Intr_Eq_1}, \eqref{Intr_Eq_2}
 (in the case  under consideration it
consists of finite or infinite curves and we shall call it graph) if all the points $\lambda \in T$
and only they
are the interior points of $G$  which are the accumulation points of
the spectrum of the problem as $k\to\infty$.
\end{Definition}

\begin{Theorem}
\label{Def_Th_2}
If both boundary points $a$ and $b$
 are finite, then the limit spectral set of
the problem \eqref{Intr_Eq_1}, \eqref{Intr_Eq_2}
 coincides with the set of balanced, essential critical
and essential singular curves, namely,
$$
T = \bigcup_{j} \gamma_j^b(a) \ \bigcup_{j} \gamma_j^a(b)\ \bigcup_{j,l}
\gamma_{j,l}^{a,b}\ \bigcup\gamma(a,b).
$$
If the point $a$ is finite but the point $b$  is infinite, then
$$
T =  \bigcup_{j} \gamma_j^b(a) \  \bigcup_{j,l} \gamma_{j,l}^{a,b}.
$$
If both points $a$ and $b$  are infinite, then:
$$
T= \bigcup_{j,l} \gamma_{j,l}^{a,b}.
$$
Here we assume that the set  $T$ also includes the endpoints of the curves which form $T$.
\end{Theorem}

The following theorem gives us formulas for the eigenvalues distribution of the problem along
the curves which constitute the limit spectral graph  $T$.

\begin{Theorem}
\label{Def_Th_3}
Let $a$ and $b$ be finite points, and the set $\gamma(a,b)$
 comprises a curvilinear interval $\gamma\subset \gamma(a,b)$
  such that
the closure of $\gamma$ has no intersections with $\gamma_j^b(a)$, $\gamma_j^a(b)$ and
 $\gamma_{,l}j^{a,b}$ Then in a neighborhood of $\gamma$ the following
quantization formula is valid
\begin{equation}
\label{Quant_1}
 k\int\limits_a^b\sqrt{P(\zeta,\lambda_m)}\,d\zeta\sim m\pi i, \quad
k\to\infty, \ m\in \mathbb Z.
\end{equation}

This  means that in a small neighborhood of the curve $\gamma$  the
eigenvalues  $\lambda_m$  of the problem
``almost'' coincide with the solutions $\lambda^0_m\in \gamma$ of the equation
$$
 k\int\limits_a^b\sqrt{P(\zeta,\lambda_m^0)}\,d\zeta= m \pi i, \quad
k\to\infty, \ m\in \mathbb Z,
$$
that is, $|\lambda_m -\lambda_m^0| = O(k^{-2})$ as $k\to\infty$.

If the point $a$  is finite, $b$ is arbitrary (finite or infinite) and there exists
a curvilinear interval $\gamma\subset \gamma_j^b(a)$
such that the  closure  of $\gamma$
does not have intersections with other curves of the graph $T$,
 then in a neighborhood of $\gamma$ the following quantization formula is valid
up to the accuracy $O(k^{-2})$:
\begin{equation}
\label{Quant_2}
 k\int\limits_{z_j(\lambda_m)}^a\sqrt{P(\zeta,\lambda_m)}\,d\zeta\sim
\left(-\frac 14 +
 m\right)\pi i,\qquad  k\to\infty,\ m\in \mathbb Z.  \notag
\end{equation}

In the case of arbitrary boundary points $a$ and $b$
 (finite or infinite),  the existence  of an curvilinear interval $\gamma\subset \gamma_{j,l}^{a,b}$
 whose  closure does not intersect the other curves of the graph $T$,  implies that  in a small
  neighborhood of $\gamma$ the following  quantization formula  is valid  up to the accuracy $O(k^{-2})$:
\begin{equation*}
\label{Quant_3}
k\int\limits_{z_j(\lambda_m)}^{z_l(\lambda_m)}\sqrt{P(\zeta,\lambda_m)}\,d\zeta\sim
\left(-\frac 12 +m\right)\pi i, \qquad k\to\infty, \ m\in \mathbb Z.
\end{equation*}
\end{Theorem}

%%%%%%%%%%%%%%%%%%%%%%%%%%%%%%%%%

As was already noted, some of the above--defined curves which form the
limit spectral graph  $T$ may also have  continual intersections. Depending on
the topology of the  Stokes graphs $\mathfrak{G}(\lambda)$ when $\lambda$ varies
 along the intervals of such intersections we
either obtain other quantization formulas (for the cases of three--points and
more intricate compound Stokes complexes) or the eigenvalues can be divided  into series
according  the
number of complexes with respect to which $a$ and $b$ are not linked. Wherein
the quantization formulas for each series are determined by the corresponding Stokes
complex. We do not present here exact formulas for such cases as we are planning
to do this in a more advanced version of the work. Also, note that the mentioned
in the theorem $O(k^{-2})$ accuracy for quantization formulas may decrease in small
neighborhoods of some critical points (here we also omit details).

At the end, we present the main ideas of the proof for
the case when the both boundary points $a$ and $b$
 are finite. First, note that all points $\lambda \in G\setminus T$
  cannot be the limit points of the spectrum when $k\to \infty$.
Indeed, if $\lambda_0\notin T$, then it follows from  the definitions  that both points  $a$ and $b$
lie in a canonical domain  $D=D(\lambda_0)$
of the complex $z$-plane (this domain is determined by the  Stokes graph
 $\mathfrak G(\lambda_0).$
In the domain $D$
 there exist  solutions $v_+$ and $v_-$
  for the equation \eqref{Intr_Eq_1} with the following asymptotic
presentations:
\begin{equation}
\label{Solut}
  v_\pm (z,\lambda,k)=\\
  \frac{1}{P^{1/4}(z,\lambda)}e^{\pm kS(z_0,z;\lambda)}(1 + O(k^{-1})), \quad \text{as} \ \, k\to\infty,
\end{equation}
uniformly with respect to $z \in K$ for any compact  $K\subset D$.

Let us take  a  domain $d$  which is compactly embedded in $D$ and contains
both points $a$ and $b$. There exists a neighborhood $U(\lambda_0)$ of the point $\lambda_0$
such that $\forall \lambda\in U(\lambda_0)$
the lines of the Stokes graphs $\mathfrak G(\lambda)$
 do not cross the domain $d$
  and the  asymptotics \eqref{Solut}
   are  valid uniformly
for all $(\lambda, z) \in U(\lambda_0) \times d$.

The characteristic determinant of the  problem \eqref{Intr_Eq_1}, \eqref{Intr_Eq_2}
has the following form:
\begin{gather}
\Delta(\lambda, k) =
\left|
\begin{array}{cc}
v_+(a) & v_+(b)\\ v_-(a) & v_-(b)
\end{array}
\right|=
P^{-1/4} (a,\lambda) P^{-1/4} (b,\lambda)\times\notag\\
\times\Bigl(
e^{kB(a,b,\lambda)}(1+O(k^{-1}))-
e^{-kB(a,b,\lambda)}(1+O(k^{-1}))\Bigr), \notag
\end{gather}
where the function $B$ is defined in \eqref{B}.
If $\lambda_0\notin \gamma(a,b)$, then $|\Delta (\lambda_0, k)|$ grows exponentially as
 $k\to \infty$; this assertion remains true
within a small neighborhood $U(\lambda_0)$, that is, $\lambda_0$ is
not a limit point of the spectrum. For $\lambda_0 \in \gamma (a,b)$, if $\lambda_0$ does  not belong to essential critical and
essential singular curves, by using the Rouche's theorem in a neighborhood of $\lambda_0$ we can get the quantization formula \eqref{Quant_1}.

In order to get quantization formulas along essential critical and singular curves one should use the transition formulas for
asymptotic solutions from one canonical domain to another one (these formulas are implemented by special transition matrices).
For example, if $\lambda_0\in \gamma_{j,l}^{a,b}$, then the points $a$ and $b$
 are located in different canonical domains which
location is determined by  two--point Stokes complexes $\Gamma_j = \Gamma_l$.
In this case the points $a$ and $b$
are located
in different basic domains of this complex, and these domains have no common boundary. Linking solutions in these
domains should be carried out by three main transition matrices:
\begin{gather}
 \begin{pmatrix} u_+(z)\\ u_-(z)\end{pmatrix} =
 \begin{pmatrix} 0 & 1\\ 1+O(k^{-1}) & i+O(k^{-1}) \end{pmatrix}
 \begin{pmatrix} 0 & \exp (k S_{j,l}(\lambda))\\ \exp (-k S_{j,l}(\lambda)) & 0
\end{pmatrix} \times\notag\\
 \times
 \begin{pmatrix} 0 & 1\\ 1+O(k^{-1}) & i+O(k^{-1})\end{pmatrix}\begin{pmatrix}
v_+(z)\\ v_-(z)\end{pmatrix},\notag
\end{gather}
where the functions $S_{jl}$ are defined in \eqref{Sjl}.
Here $u_+(z)$ and $u_-(z)$ form a fundamental system of solutions in the basic domain containing
 point $b$, while $v_+(z)$ and $v_-(z)$ form a fundamental system of solutions with asymptotics  \eqref{Solut}
 in the basic domain  containing the point  $a$. Using this relation for
calculating  the characteristic determinant we get \eqref{Quant_3}. Similarly (or even
more simpler), the formulas \eqref{Quant_2} can be obtained.

\bigskip
{\Large\bf\noindent Examples}
\bigskip

To illustrate the results obtained above let us consider two examples.

{\bf Example 1.} (cf. \cite{Sh1}). {\it The model problem for plane Couette flow.}
\begin{equation}
\label{example1}
-y''+k^2i(z-\lambda)y=0,\quad y(-1)=y(1)=0.
\end{equation}

In this case $P(z,\lambda)=i(z-\lambda)$. For each $\lambda\in\CC$ there is the only 
 single turning point $z_1=\lambda$. One can explicitly calculate
$$
S(z_1,z;\lambda)=\int\limits_{\lambda}^z\sqrt{i(\zeta-\lambda)}\,d\zeta=
\frac{2}{3}e^{i\pi/4}(z-\lambda)^{3/2}.
$$

\begin{figure}
\begin{center}
\includegraphics[width=15.9cm,keepaspectratio]{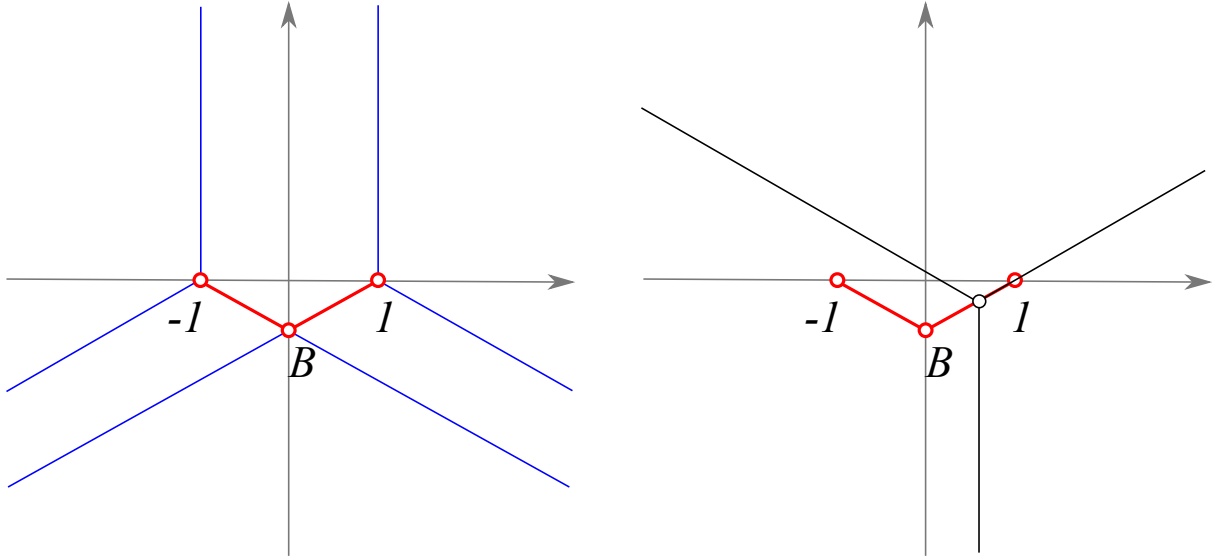}
\end{center}
\caption{Example 1. On the left: Critical (blue) and essential critical (red) curves; on the right: Stokes Graph (black) for $\lambda$ lying on
the essential critical curve.}
\label{graph1}
\end{figure}

For each $\lambda$ the Stokes Graph consists of a single complex whith Stokes lines which are straight-line rays starting from $\lambda$:
$$
\Gamma(\lambda)=\Bigl\{
\zeta=\lambda+re^{i\psi_k}\,\Bigl|\Bigr.\,
r\ge0,\ \psi_k=\frac{\pi}{6}+\frac{2\pi k}{3},\ k=0,1,2
\Bigr\}.
$$
From this, in particular, we get the absence of essential singular curves in the limit spectral set. The location of critical and
essential critical curves is illustrated on fig.\ref{graph1}. Here $B=-i/\sqrt{3}$ is the intersection point of the critical curves.

Basing on the nature of balanced curves we get that $\pm1$ should not only belong to common canonical domain, but furthermore to one
basic domain which in our case is valid only in domain $D$, which lies below all critical curves:
$$
D=\Bigl\{
\lambda\in\CC
\,\Bigl|\Bigr.\,
\lambda=-i/\sqrt{3}+\xi,\ \arg\xi\in(-\frac{5\pi}{6};-\frac{\pi}{6})
\Bigr\},
$$
and the equation for balanced curves has the form: $\re S(\lambda,1;\lambda)=\re S(\lambda,-1;\lambda)$, $\lambda\in D$. Its solution is
the infinite interval $\gamma_\infty=(-i/\sqrt{3};-i\infty)$.

\begin{figure}
\begin{center}
\includegraphics[width=15.9cm,keepaspectratio]{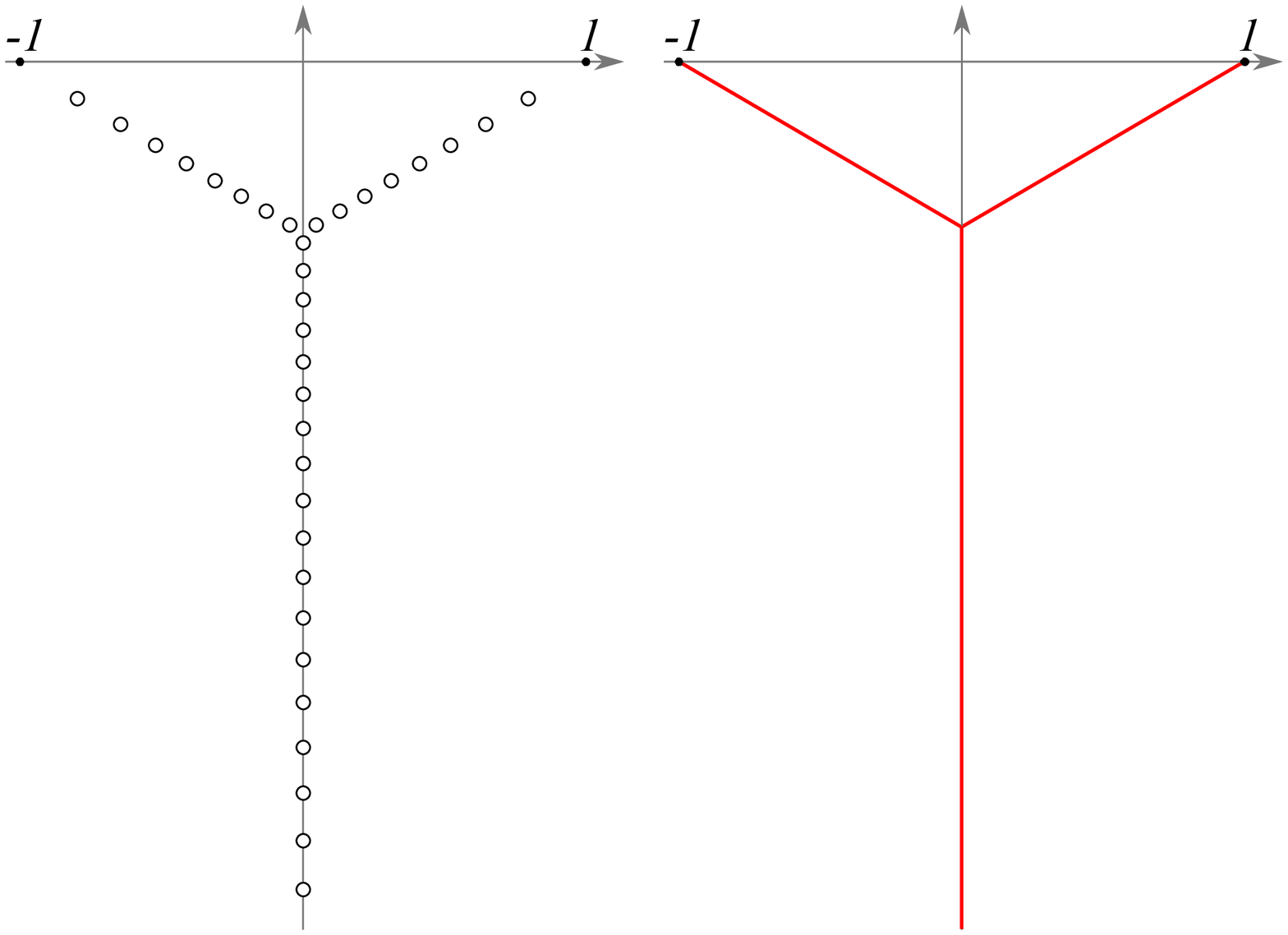}
\end{center}
\caption{Example 1. On the left: Eigenvalues for some $k>0$. On the right: Limit spectral graph (red).}
\label{graph2}
\end{figure}

Eigenvalues of \eqref{example1} for some sufficiently large $k$ and the limit spectral set are shown on fig.\ref{graph2}.

\bigskip
{\bf Example 2.} {\it The model problem for plane Couette--Poiseuille flow.}

Consider the problem:
\begin{equation}
\label{example2}
-y''+k^2i(z^2-\frac{1}{2}z-\lambda)y=0,\quad y(-1)=y(1)=0.
\end{equation}

Now $P(z,\lambda)=i(z^2-z/2-\lambda)$. Both turning points $z_{1,2}=1/4\pm\sqrt{1/4+4\lambda}$ are simple when $\lambda\ne-1/16$.

The equation for singular curves is as follows:
$$
\re\int\limits_{z_1}^{z_2}\bigl(i(\zeta-z_1)(\zeta-z_2)\bigr)^{1/2}\,d\zeta=0.
$$
 It is equivalent to the equation  $\im e^{\pi i/4}(z_2-z_1)^2=0$, i.e. $\im e^{\pi i/4}(1/4+4\lambda)=0$. The singular curve is a staight line in $\CC$:
$\im\lambda+\re\lambda+1/16=0$.

Let us multiply the basic equation \eqref{example2} by $\overline{y}$ and integrate the obtained equation 
from -1 to 1. Then, separating the real and imaginary parts, we obtain:
$$
\im\lambda=-\frac{1}{k^2}\frac{\langle y',y'\rangle}{\langle y,y\rangle}\subset(-\infty,0),\quad
\re\lambda=\frac{\langle(x^2-x/2)y,y\rangle}{\langle y,y\rangle}\subset(-1/16,3/2).
$$

 This analysis allows us to narrow the domain $G$ where we study the spectrum. In our case it is  reduced 
  from the whole plane to the half--strip:
$$
G=\Bigl\{
\lambda\in\CC
\,\Bigl|\Bigr.\,\re\lambda\in(-1/16,3/2),\ \im\lambda<0
\Bigr\}.
$$

\begin{figure}
\begin{center}
\includegraphics[width=15.9cm,keepaspectratio]{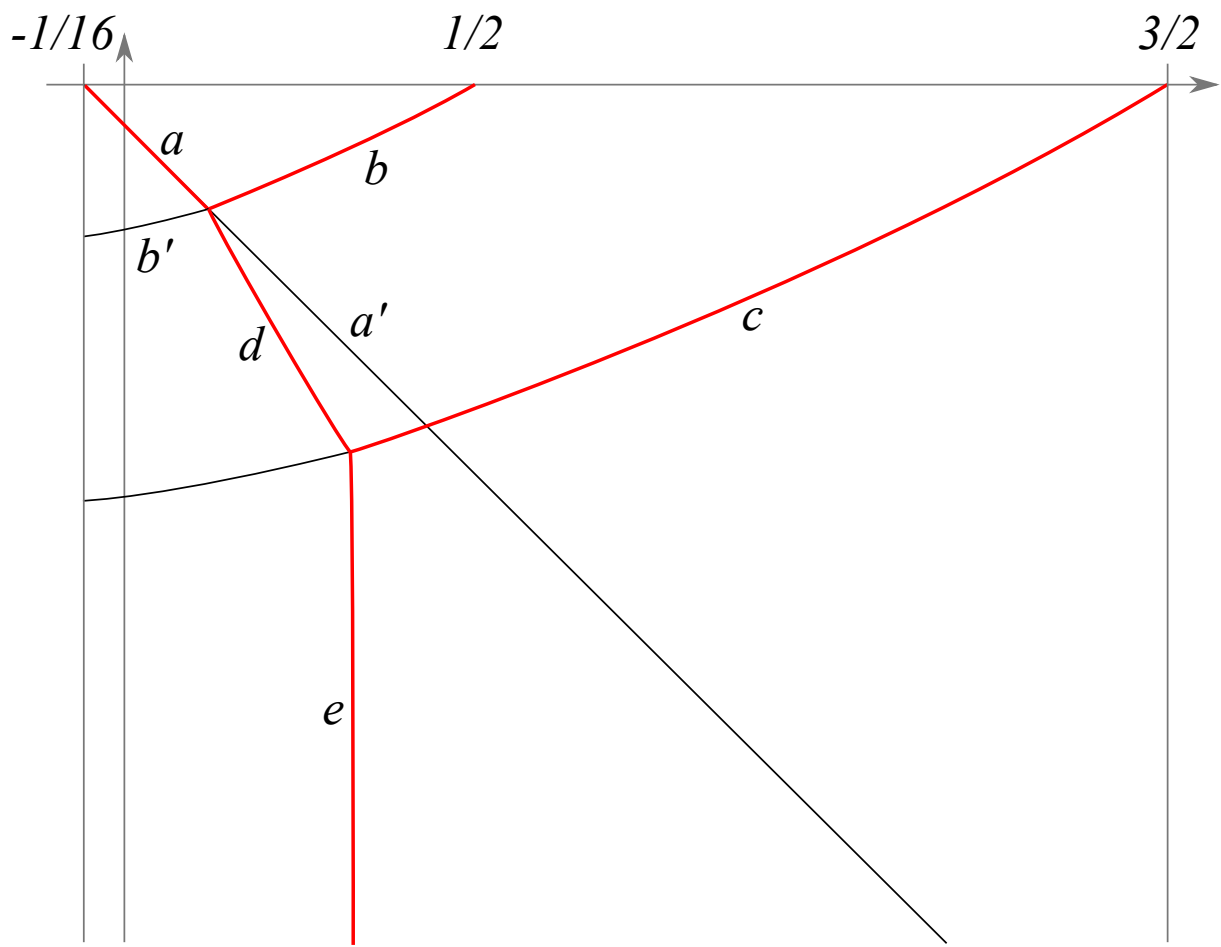}
\end{center}
\caption{Example 2. Limit spectral set (red); examples of singular and critical curves which are not essential singular and not essential critical (black).}
\label{graph3}
\end{figure}

\begin{figure}
\begin{center}
\includegraphics[width=15.9cm,keepaspectratio]{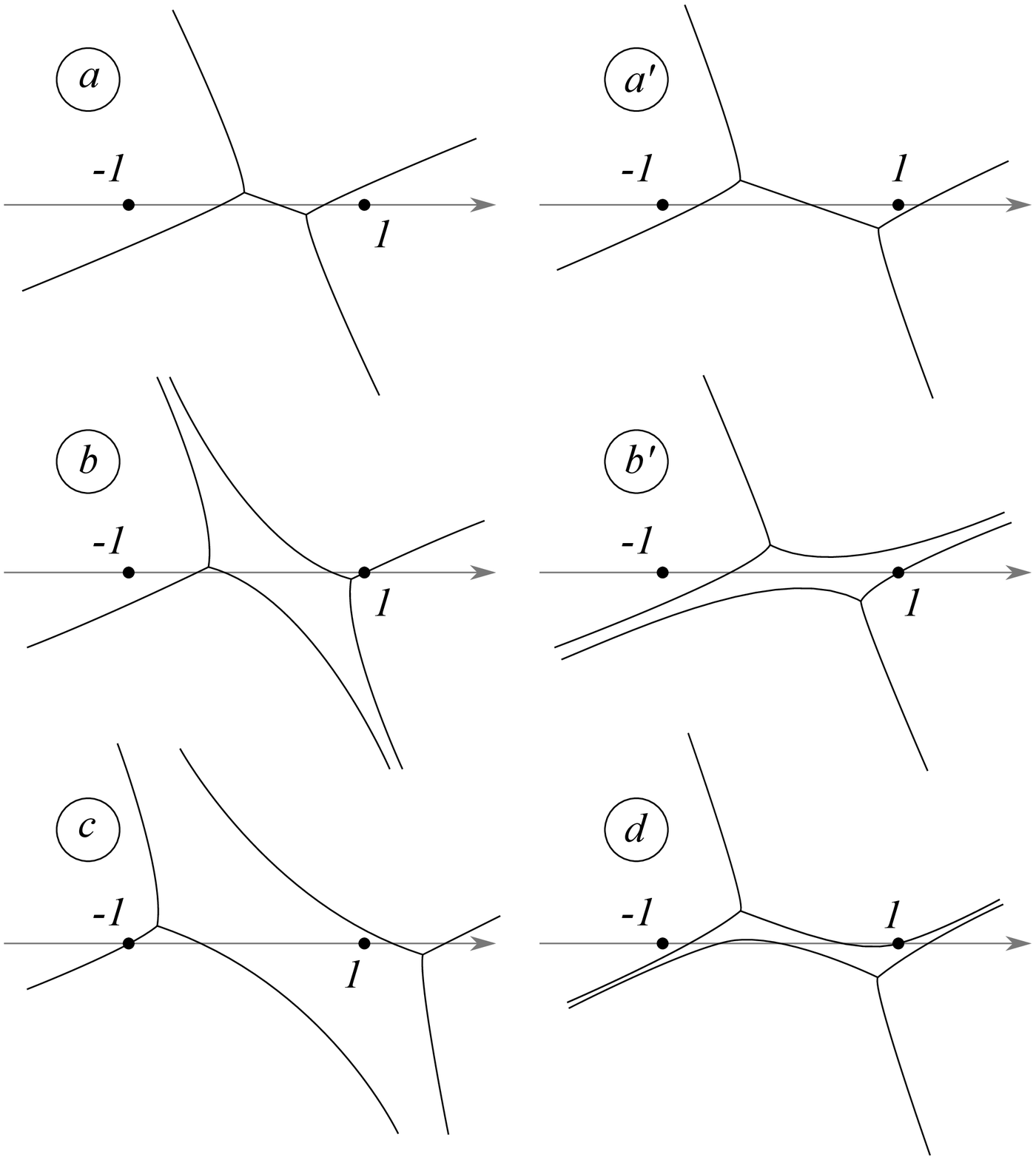}
\end{center}
\caption{Example 2. Location of $\pm1$ with respect to to Stokes Graph.}
\label{graph4}
\end{figure}

Figure \ref{graph3} illustrates the limit spectral set (in red): essential singular curve $(a)$, essential critical curves:
$(b)$, $(c)$, $(d)$ and the balanced curve $(e)$. We also put there an example of singular, but not essential singular curve
$(a')$ and critical, but not essential critical curve $(b')$: along these curves $\pm1$ are linked with respect to all Stokes complexes,
for this reason they are not contained in the limit spectral set.

Figure \ref{graph4} illustrates the location of $\pm 1$ with respect to Stokes Graph. Evidently $\pm 1$ are not linked in the case $(a)$.
In the case $(b)$ $\pm1$ are linked with respect to the left complex, but not linked with respect to the right one. On the contrary in the case
$(c)$ $\pm1$ are not linked with respect to the left complex, but linked with respect to the right one. In the case $(d)$ $\pm1$ are not
linked with respect to the upper complex, but linked with respect to the lower one.


\begin{thebibliography}{99}

\bibitem{TSh4}
A.\,A.~Shkalikov and S.\,N.~Tumanov. {\it
The Limit Spectral Graph in Quasi-Classical Approximation
for the Sturm–Liouville Problem with Complex Polynomial Potential}
Doklady Mathematics, 2015, Vol. 92, No. 3.

\bibitem{EvgrafFedoruk}
M.\,A.~Evgrafov and M.\,V.~Fedoryuk. {\it Asymptotic behavior as $\lambda\to +\infty$ of solutions of the equation $w''(z)-p(z,\lambda)w(z)=0$
in the complex $z$--plane.} Russ. Math. Surveys, \textbf{21}, 1 (1966), 3--50.

\bibitem{ESh}
A.\,I.~Esina and A.\,I.~Shafarevich. {\it Quantization conditions on Riemannian surfaces
and the Semi--Classical Spectrum of the Schr\"odinger operator
with complex potential.} Math. Notes, {\bf 88}:2
(2010), 229--248.

\bibitem{PSh}
V.\,I.~Pokotilo and A.\,A.~Shkalikov. {\it Semiclassical approximation for a nonself--adjoint Sturm--Liouville problem with a parabolic
potential.} Math. Notes, {\bf 86}:3
(2009), 561--569.

\bibitem{TSh1}
A.\,A.~Shkalikov and S.\,N.~Tumanov. {\it On the limit behaviour of the spectrum of a model problem for the Orr--Sommerfeld equation with Poiseuille profile.}
Izv. Math., {\bf 66}:4 (2002), 177--204.

\bibitem{TSh2}
A.\,A.~Shkalikov and S.\,N.~Tumanov. {\it On the spectrum localization of the Orr--Sommerfeld problem for large Reynolds numbers.}
Math. Notes, {\bf 72}:4
(2002), 561--569.

\bibitem{TSh3}
A.\,A.~Shkalikov and S.\,N.~Tumanov. {\it On the model for the Orr--Sommerfeld equation with quadratic profile.} arXiv: math-ph/0212074v1, 2002.

\bibitem{Fedoruk1}
M.\,V.~Fedoryuk. {\it Asymptotic methods for linear ordinary differential equations.}
(Mir, Moscow, 1983) [in Russian].

\bibitem{Fedoruk2}
M.\,V.~Fedoryuk. {\it Topology of the Stokes lines of a second-order equation.} Izv. Akad. Nauk SSSR, Ser. Mat., \textbf{29}, 3 (1965),
645--656.

\bibitem{Fedoruk3}
M.\,V.~Fedoryuk. Addition I to W.~Wazow, {\it Asymptotic Expansions for Ordinary Differential
Equations.} (Mir, Moscow, 1968) [in Russian].

\bibitem{Sh1}
A.\,A.~Shkalikov. {\it The limit behavior of the spectrum for large parameter values in a model problem.} Math. Notes, {\bf 58}:6 (1997), 950--953.

\bibitem{EG}
A.~Eremenko and A.~Gabrielov. {\it Singluar perturbation of polynomial potentials in the complex domain with applications to
spectral loci of $\mathcal{PT}$--symmetric families}. Mosc. Math. J., {\bf 11}:3 (2011),
473--503.

\bibitem{Heding}
J.\,Heading. {\it An introduction to phase--integral methods.}
(Wiley, New York, 1962; Mir, Moscow, 1965).

\bibitem{Sh2}
A.\,A.~Shkalikov. {\it Spectral portraits of the Orr--Sommerfeld operator
with large reynolds numbers.} J. Math Sci., {\bf 124}:6 (2004),
5417--5441.
\end{thebibliography}
\end{document}